\documentclass[11pt]{amsart} \usepackage{amsfonts,amsmath,amssymb}

\newtheorem{definition}{Definition}[section]

\newtheorem{proposition}[definition]{Proposition}
\newtheorem{corollary}[definition]{Corollary}

\newtheorem{theorem}[definition]{Theorem}
\newtheorem{example}[definition]{Example}
\def\bea{\begin{eqnarray*}}\def\eea{\end{eqnarray*}}

\newcommand{\Hom}{{\rm Hom}}
\newcommand{\End}{{\rm End}}
\newcommand{\can}{{\rm can}}

\def\a{$\mbox{\u{a}}$}

\def\S{$\mbox{\c{S}}$}

\def\bea{\begin{eqnarray*}}\def\eea{\end{eqnarray*}}

\begin{document}
\title{Duality for finite Hopf algebras explained by corings} \author{S. Caenepeel}\thanks{This paper was written while the first author visited the Mathematics Departments of Syracuse University and California State University Dominguez Hills. He would like to thank both departments for their hospitality.} \address{Vrije Universtiteit Brussel (VUB), Pleinlaan 2, B-1050, Brussel, 
Belgium}
\email{scaenepe@vub.ac.be}
\author{D. Quinn}
\address{Department of Mathematics,
Syracuse University, Syracuse, NY 13244, USA} \email{dpquinn@syr.edu} \author{\c{S}. Raianu} \address{Mathematics Department, California State University, Dominguez Hills, 1000 E Victoria St, Carson CA 90747, USA} \email{sraianu@csudh.edu} 
\date{May 4, 2005}

\begin{abstract}
We give a coring version for the duality theorem for actions and coactions of a finitely generated projective Hopf algebra. We also provide a coring 
analogue for a theorem of 
H.-J. Schneider, which generalizes and unifies the duality theorem for finite Hopf algebras and its refinements. \end{abstract} \maketitle \section*{Introduction} The Duality Theorem for actions and coactions of a finite Hopf algebra is the following: \begin{theorem}\label{dual}\cite{bm, vdb} Let $H$ be a finitely generated projective Hopf algebra over the commutative ring $k$, and $A$ a right $H$-comodule algebra. Then $$(A\# H^*)\# H\simeq \End_A(H\otimes A).$$ \end{theorem} The following are two refinements of Theorem \ref{dual}: \begin{theorem}\label{ulbrich}\cite{u}
Let $H$ be a finitely generated projective Hopf algebra over the commutative ring $k$, and $A$ a right $H$-comodule algebra. Then:\\
i) $\End^H_A(H\otimes A)\# H\simeq \End_A(H\otimes A)$\\
ii) $\End^H_A(H\otimes A)\simeq A\#H^*.$
\end{theorem}
It was shown in \cite{drv} that both assertions of the previous theorem are particular cases of the following result, due to H.-J. Schneider: \begin{theorem}\label{hans}\cite{sch}
Let $H$ be a finitely generated projective Hopf algebra over the commutative ring $k$, $A$ a left $H$-module algebra, and $B=A^H$ the algebra of invariants. Assume that $A/B$ is a right $H^*$-Hopf Galois extension, and let $M\in {\mathcal M}^{H^*}_A$. Then $End_A(M)$ is a left $H$-module algebra via the action $(h\cdot f)(m)= h_{(1)}f(S(h_{(2)})m)$, and we have the algebra isomorphism $$\End_A(M)\# H\simeq \End_B(M),$$ sending $f\# h$ to the endomorphism of $M$ mapping $m$ to $f(hm)$. 
\end{theorem} Finding the point of view that makes the proof of the finite duality theorem trivial has at least two possible answers: one could be ``Hopf-Galois extensions'', as in \cite{vdb}, and the second one could be ``the fundamental theorem for Hopf modules'', because this result was used in \cite{drv} to prove a more general form of Schneider's theorem.

There is yet another question about the finite duality theorem that seems worth answering. There are apparently two distinct types of duality results.

The first one (like in Theorem \ref{dual}) could be described as follows: take an algebra; let $H^*$ act on it and take the smash product; $H$ acts on it, so take the smash product again; get back to the algebra (actually to its endomorphism ring, which is Morita equivalent to it).

The second one is: take object; take dual; take dual again; get back to object.

By comparing the two, we see that the way to make them agree is to describe the smash product as the dual of an object. This object is a coring. \begin{definition} Let $A$ be a ring. An $A$-bimodule ${\mathcal C}$ is called a {\em coring} if there exist $A$-bimodule maps $\Delta:{\mathcal C}\longrightarrow {\mathcal C}\otimes_A{\mathcal C}$ and $\varepsilon:{\mathcal C}\longrightarrow A$, such that $\Delta$ is coassociative and $\varepsilon$ is a counit.  
\end{definition}
Corings were introduced by Sweedler in \cite{swe}, and were given a lot of attention lately, after Takeuchi remarked that many examples of (generalized) Hopf modules are in fact just comodules over some corings. For all unexplained facts about corings the reader is referred to \cite{BW}. 

The most basic example of a coring is a coalgebra over a commutative ring $A$. Another fundamental example is the canonical coring associated to the ring homomorphism $i:B\longrightarrow A$: ${\mathcal C}=A\otimes_B A$, $\Delta:{\mathcal C}\longrightarrow{\mathcal C}\otimes_A{\mathcal C}\simeq A\otimes_B A\otimes_B A$, $\Delta(a\otimes b)=a\otimes 1\otimes b$, and $\varepsilon:{\mathcal C}\longrightarrow A$, $\varepsilon(a\otimes b)=ab$.

\begin{example}\label{coringatensorh}
Let $H$ be a $k$-Hopf algebra, and $A$ a right $H$-comodule algebra via $a\mapsto a_{[0]}\otimes a_{[1]}$. Then $A\otimes H$ becomes an $A$-coring as
follows: the left $A$-module structure is given by multiplication on the first component, and the right $A$ module structure is given by $(a\otimes h)b= ab_{[0]}\otimes hb_{[1]}$. The comultiplication is $\Delta:A\otimes H \longrightarrow (A\otimes H)\otimes_A(A\otimes H)\simeq A\otimes H\otimes H$, $\Delta(a\otimes h)=a\otimes h_{(1)}\otimes h_{(2)}$, and the counit is $\varepsilon:A\otimes H\longrightarrow A$, $\varepsilon(a\otimes h)= \varepsilon(h)a.$ \end{example}

Recall that if $A$ is a right $H$-comodule algebra, then the smash product $A\# H^*$ has multiplication defined by $(a\# h^*)(b\# g^*)=ab_{[0]}\# (h^*\leftharpoonup b_{[1]})g^*$. The big smash product $\#(H,A)=\Hom_k(H,A)$ has multiplication defined by $(f\cdot g)(h)=f(g(h_{(2})_{[1]}h_{(1)}) g(h_{(2})_{[0]}$. If $H$ is projective finitely generated over $k$, and $h_i$ and $h^*_i$, $i=1\ldots n$, are dual bases, then $\#(H,A)\simeq A\# H^*$ as rings, via $f\mapsto \sum_{i=1}^n f(h_i)\# h^*_i$. Note that if $B=A^{coH}$, then $A/B$ is a right $H^*$-Hopf-Galois extension if and only if the coring $A\otimes H$ from example \ref{coringatensorh} is isomorphic to the canonical coring $A\otimes_B A$ via the map sending $a\otimes b\in A\otimes_B A$ to $ab_{[0]}\otimes b_{[1]}\in A\otimes H$.

\section{Duality for corings}
In this section we prove a duality theorem for finitely generated projective  corings, and show that it extends the duality theorem for finite Hopf algebras.\\ Let ${\mathcal C}$ be an $A$-coring. Then the left dual of ${\mathcal C}$, 
denoted by ${^*{\mathcal C}}$, is by definition
$${^*{\mathcal C}}={_A{\Hom({\mathcal C},A)}},$$
with multiplication given by $$(\phi*\psi)(c)=\psi(c_{(1)}\phi(c_{(2)})).$$
Then $A$ becomes a subring of ${^*{\mathcal C}}$ via $$i:A\longrightarrow {^*{\mathcal C}},\;\; i(a)(c)=\varepsilon(c)a.$$ Similarly, the right dual of ${\mathcal C}$, 
denoted by ${\mathcal C}^*$, is by definition
$${\mathcal C}^*=\Hom_A({\mathcal C},A),$$
with multiplication given by $$(\phi*\psi)(c)=\phi(\psi(c_{(1)})c_{(2)}).$$
Now $A$ becomes a subring of ${\mathcal C}^*$ via $$j:A\longrightarrow {\mathcal C}^*,\;\; j(a)(c)=a\varepsilon(c).$$

\begin{theorem}\label{dualcoring}
Let ${\mathcal C}$ be an $A$-coring, and assume that ${\mathcal C}$ 
is finitely generated projective as a left $A$-module. Then the right dual of the canonical coring associated to the ring morphism $i:A\longrightarrow {^*{\mathcal C}}$ is isomorphic to the (opposite) endomorphism ring 
$\End({_A{\mathcal C}})^{op}$.
\end{theorem}

\begin{proof}
We have the isomorphism
$$\alpha:\ \Hom_{^*{\mathcal C}}({^*{\mathcal C}}\otimes_A{^*{\mathcal C}}, {^*{\mathcal C}}) \stackrel{\textstyle{\longrightarrow}}{\longleftarrow}
\Hom_A({^*{\mathcal C}},{^*{\mathcal C}}):\beta,$$
defined by $\alpha(\phi)(c^*)=\phi(c^*\otimes 1)$ and $\beta(\psi)(c^*\otimes d^*)=\psi(c^*)d^*$. (Note that $\alpha$ is simply the compostion of the hom-tensor adjointness with the equivalence of $Hom_R(R,-)$, for $R = ^*{\mathcal C}$ 
an associative ring.) Now duality for finitely generated projective modules provides the isomorphism $$\gamma:\ \Hom_A({^*{\mathcal C}},{^*{\mathcal C}}) \stackrel{\textstyle{\longrightarrow}}{\longleftarrow}
{_A{\Hom({\mathcal C},{\mathcal C})}}:\delta,$$
which is defined as follows: let $c=\sum_{i=1}^nf_i(c)c_i$ be the dual basis formula in ${\mathcal C}$; then $\gamma(\psi)(c)=\sum_{i=1}^n\psi(f_i)(c)c_i$,
and $\delta(\zeta)(f)(c)=f(\zeta(c))$.
It is easy to check that $\alpha$ is a ring isomorphism, while $\delta$ is a ring anti-isomorphism. \end{proof}

Let us show now how the duality theorem for finitely generated projective Hopf algebras can be derived from Theorem \ref{dualcoring}. We first need a 
description of the smash product as a dual of a coring. Before giving this result, let us remark that a coring interpretation of Theorem \ref{ulbrich} ii) also realizes the smash product as the endomorphism ring of a coring. 
\begin{proposition}\label{smashdualcoring}
Let $H$ be a finitely generated projective Hopf alegbra over the commutative ring $k$, and $A$ a right $H$-comodule algebra. Let $A\otimes H$ be the 
$A$-coring from example \ref{coringatensorh}. There there are ring isomorphisms $${^*{(A\otimes H)}}\simeq(A\otimes H)^*\simeq \#(H,A)\simeq A\# H^*.$$ \end{proposition}
\begin{proof}
Define
$$\gamma:(A\otimes H)^* \stackrel{\textstyle{\longrightarrow}}{\longleftarrow}
\#(H,A):\delta$$
by $\gamma(\phi)(h)=\phi(1\otimes S^{-1}(h))$, and $\delta(\zeta)(a\otimes h)= \zeta(a_{[1]}S(h))a_{[0]}$. Since it is easy to see that $\gamma(\phi)$ is $k$-linear and $\delta(\zeta)$ is right $A$-linear, and that they are inverse one to each other, we check that they are ring isomorphisms: \begin{eqnarray*} &&\hspace*{-2cm} (\gamma(\phi)\cdot\gamma(\psi))(h)=
\gamma(\phi)(\gamma(\psi)(h_{(2)})_{[1]}h_{(1)})\gamma(\psi)(h_{(2)})_{[0]}\\
&=&\phi(1\otimes S^{-1}(\psi(1\otimes S^{-1}(h_{(2)})_{[1]}h_{(1)}) \psi(1\otimes S^{-1}(h_{(2)})_{[0]}\\ &=&\phi(\psi(1\otimes S^{-1}(h_{(2)})\otimes S^{-1}(h_{(1)})\\ &=&\phi(\psi(1\otimes S^{-1}(h_{(2)})(1\otimes S^{-1}(h_{(1)}))\\ &=&\gamma(\phi*\psi)(h). \end{eqnarray*} This completes the proof for the middle isomorphism. Now define $$\gamma':{^*{(A\otimes H)}} \stackrel{\textstyle{\longrightarrow}}{\longleftarrow}
\#(H,A):\delta'$$
by $\gamma'(\phi)(h)=\phi(1\otimes h)$, and $\delta'(\zeta)(a\otimes h)= a\zeta(h)$. It is clear that $\gamma'(\phi)$ is $k$-linear and $\delta'(\zeta)$ is left $A$-linear, and that they are inverse one to each other, so we check that they are ring isomorphisms. A short computation shows that $\gamma'$ establishes a ring isomorphism ${^*{(A\otimes H)}}\simeq \#(H^{op},A^{op})^{op}\simeq 
(A^{op}\# H^{*cop})^{op}$. Finally, 
$$A\# H^*\simeq (A^{op}\# H^{*cop})^{op}, 
\;\; a\# h^*\mapsto h^*_{(1)}(S^{-1}(a_{[1]}))a_{[0]}\# h^*_{(2)}S^{-1}$$ 
as rings, and the proof is complete. 
\end{proof}

\begin{corollary}\label{dualityhopf}
Let $H$ be a finitely generated projective Hopf algebra over $k$, and $A$ a 
right $H$-comodule algebra. Then $(A\# H^*)\# H\simeq \End({A\# H^*}_A)$ as rings. \end{corollary}

\begin{proof}
We have the following diagram of isomorphisms:

\begin{picture}(100,110)(0,0) \put(50,10){$\Hom_k(H^*,{^*{\mathcal C}})$} \put(50,45){$\Hom_{^*{\mathcal C}}({^*{\mathcal C}}\otimes H^*, {^*{\mathcal C}})$} \put(50,80){$\Hom_{^*{\mathcal C}}({^*{\mathcal C}}\otimes_A {^*{\mathcal C}}, {^*{\mathcal C}})$} \put(200,80){$\Hom_A({^*{\mathcal C}}, {^*{\mathcal C}})$} \put(150,85){\vector(1,0){35}} \put(75,22){\vector(0,1){20}} \put(75,55){\vector(0,1){20}} \put(130,15){\vector(4,3){65}} \end{picture}

\noindent where ${\mathcal C}$ is the $A$-coring $A\otimes H$ from Example 
\ref{coringatensorh}, and 
the diagonal map is the composition of the two vertical maps and the 
horizontal one. 
By Proposition \ref{smashdualcoring}, the lower corner
of the diagram is just $(A\# H^*)\# H$. Also by Proposition 
\ref{smashdualcoring}, ${^*{\mathcal C}}$ is a right $H^*$-comodule algebra. 
Let us denote this structure by $c^*\mapsto c^*_{[0]}\otimes c^*_{[1]}$. 
The horizontal isomorphism is the one from
the first part of the proof of Theorem \ref{dualcoring}, the lower vertical one assigns to $f\in \Hom_k(H^*,{^*{\mathcal C}})$ the map sending $c^*\otimes h^*$ to $f(c^*_{[1]}h^*S)c^*_{[0]}$ (like $\zeta$ from Proposition \ref{smashdualcoring}). The upper vertical isomorphism is induced by the map $$ \can:\ {^*{\mathcal C}}\otimes_A {^*{\mathcal C}}\longrightarrow {^*{\mathcal C}}\otimes H^*,\;\; c^*\otimes d^*\mapsto c^*d^*_{[0]}\otimes d^*_{[1]}.$$ Thus the diagonal map in the above diagram provides the isomorphism in the 
statement. In order to make sure that this is exactly the isomorphism from the duality theorem for finite Hopf algebras, we only need to show that it assigns to $f$ the map sending $c^*$ to $f(c^*_{[1]})c^*_{[0]}$, which checks immediately. \end{proof}

In view of the preceding Corollary, we get yet another possible point of view that trivializes the duality theorem for finite Hopf algebras: behind it we find just duality for finitely generated projective modules.

We also have an explanation for getting the bigger endomorphism ring of $A$, rather than $A$ itself after taking the dual for the second time. Once we take the dual of a coring, we get a ring, so to get another coring that will enable us to take the second dual we need another construction, namely passing to the canonical coring. 

\section{A coring version of Schneider's isomorphism}
Recall that if $({\mathcal C},x)$ is an $A$-coring with fixed grouplike element $x$, and $B=A^{co{\mathcal C}}=\{a\in A\mid \Delta(a)=a\otimes x\}$, then $({\mathcal C},x)$ is a Galois coring if the canonical coring morphism $can:A\otimes_B A\longrightarrow {\mathcal C}$, $\can(a\otimes b)=axb$ is an isomorphism. We have the following coring version of Theorem \ref{hans}:

\begin{theorem}\label{hansgen}
Let $({\mathcal C},x)$ be a Galois coring, $M\in{\mathcal M}^{\mathcal C}$, and $N\in {\mathcal M}_A$. Then 
$$\Hom_A(M\otimes _A{\mathcal C},N)\simeq \Hom_B(M,N).$$ \end{theorem}

\begin{proof}
\begin{eqnarray*}
&&\hspace*{-2cm}
\Hom_A(M\otimes _A{\mathcal C},N)  \simeq \Hom_A(M\otimes _AA\otimes_BA,N)\\ & \simeq & \Hom_A(M\otimes_BA,N)\\ & \simeq & \Hom_B(M,\Hom_A(A,N))\\ 
&&(\mbox{by the adjunction property of the tensor product)}\\
& \simeq & \Hom_B(M,N).\\
\end{eqnarray*}
\end{proof}

Theorem \ref{hans} follows now immediately from Theorem \ref{hansgen}: under the hypotheses of Theorem \ref{hans} (with $H$ replaced by $H^*$), take ${\mathcal C}=A\otimes H$ as in Example \ref{coringatensorh}, $x=1\otimes 1$, and $M=N\in{\mathcal M}^{\mathcal C}$. Then we have that $H\otimes M\simeq M\otimes H$ as right $A$-modules, and therefore:
\begin{eqnarray*}
End_B(M) & \simeq & Hom_A(M\otimes _A(A\otimes H),M)\\
 & \simeq & Hom_A(M\otimes H,M)\\
& \simeq & Hom_A(H\otimes M,M)\\
& \simeq & Hom_k(H,Hom_A(M,M))\\ 
&&(\mbox{by the adjunction property of the tensor product)}\\
& \simeq & End_A(M)\otimes H^*.
\end{eqnarray*}

We leave it to the reader to check that the isomorphism is the one from Theorem \ref{hans} (just note that the last isomorphism from above sends $g$ to $\sum_{i=1}^n g(h_i)\# h^*_iS^{-1}$).

In view of the above, we also obtain the surprising fact that what is behind the duality theorem for finite Hopf algebra and its refinements is just the hom-tensor adjunction property.


\thebibliography{MMM} 
\bibitem{abe} E. Abe, {\it Hopf Algebras}, Cambridge Univ. Press, 1977. 
\bibitem{bm} R.J. Blattner, S. Montgomery, A duality theorem for Hopf module algebras. J. Algebra {\bf 95} (1985), no. 1, 153-172. 
\bibitem{BW} T. Brzezi\'nski and R. Wisbauer, {\it Corings and comodules}, 
London Math. Soc. Lect. Note Ser. {\bf 309},
  Cambridge University Press, Cambridge, 2003.
entwining structures, {\sl J. Algebra} {\bf 276} (2004), 210--235. \bibitem{dnr} S. D\a sc\a lescu,  C. N\a st\a sescu, \S. Raianu, {\it Hopf Algebras: an Introduction}, Monographs and Textbooks in Pure and Applied Mathematics {\bf 235}, 
Marcel Dekker, Inc., New York, 2001.
\bibitem{drv}
S. D\a sc\a lescu, \S. Raianu, F. Van Oystaeyen,
Some Remarks on a Theorem of H.-J. Schneider, Comm. Algebra {\bf 24}(14), 
(1996), 4477--4493.
\bibitem{mr}
C. Menini, \S. Raianu, Morphisms of relative Hopf modules, smash products and duality, J. Algebra {\bf 219} (1999), 547--570. \bibitem{mont} S. Montgomery, {\it Hopf algebras and their actions on rings}, CBMS Regional Conference Series in Mathematics {\bf 82}, 
Amer. Math. Soc., Providence, RI, 1993.
\bibitem{sch}
H.-J. Schneider, Hopf Galois extensions, crossed products and Clifford Theory, in {\it Advances in Hopf Algebras}, edited by J. Bergen and S. Montgomery, Marcel Dekker Lecture Notes in Pure and Appl. Math., {\bf 158} (1994), 267--298. \bibitem{sw} 
M.E. Sweedler, {\it Hopf Algebras}, Benjamin, New York, 1969. \bibitem{swe} M.E. Sweedler, The predual theorem to the Jacobson-Bourbaki theorem, Trans. Amer. Math. Soc. {\bf 213} (1975), 391--406. \bibitem{u} K.-H. Ulbrich, Smash products and comodules of linear maps, Tsukuba J. Math., {\bf 14} (1990), 371--378. \bibitem{vdb} M. Van den Bergh, A duality theorem for Hopf algebras, in Methods in ring theory (Antwerp, 1983), 517--522, NATO Adv. Sci. Inst. Ser. C Math. Phys. Sci., {\bf 129}, Reidel, Dordrecht, 1984.
\end{document}